# Simulation of thin film flows with a moving mesh mixed finite element method

Hong Zhang*, Paul Andries Zegeling

*Department of Mathematics, Faculty of Science, Utrecht University, Budapestlaan 6, 3584CD Utrecht, The Netherlands*

**Abstract**

We present an efficient mixed finite element method to solve the fourth-order thin film flow equations using moving mesh refinement. The moving mesh strategy is based on harmonic mappings developed by Li et al. [J. Comput. Phys., 170 (2001), pp. 562-588, and 177 (2002), pp. 365-393]. To achieve a high quality mesh, we adopt an adaptive monitor function and smooth it based on a diffusive mechanism. A variety of numerical tests are performed to demonstrate the accuracy and efficiency of the method. The moving mesh refinement accurately resolves the overshoot and downshoot structures and reduces the computational cost in comparison to numerical simulations using a fixed mesh.

## 1. Introduction

Thin liquid film flows appear in various situations in nature and engineering applications, such as rain flow down along a window, spin coating, lubrication, membranes in biophysics, etc. Despite the diversity of applications, the governing model is similar if the film is sufficiently viscous. Huppert first explained the fingering instability of the thin film flows using a mathematical model in [1]. After Huppert's initial research on liquid film flow down an incline, experiments have been carried out in different configurations, such as the rewetting of an inclined solid surface [2], thermally driven (Marangoni effect) thin films [3] and so on. The experimental studies reveal that, in the case of forced spreading, the thin film front undergoes a fingering instability. Besides experimental studies, theoretical researches of thin film flows have also been conducted in the last few decades under different aspects. Results on solution existence, traveling wave, phase plane analysis, pattern formation, stability, kinetics and nucleation are given in [4, 5, 6, 7, 8, 9].

In order to study the behaviors of the thin film flow equations numerically, a variety of numerical methods have been developed in literature. Bertozzi and Bowen [10] implemented a positivity preserving finite difference scheme. Ha et al. [11] studied the stability of traveling wave solutions to the thin film equation by comparing the solutions obtained by Crank-Nicolson, fully implicit, Godunov, adapted upwind and weighted essentially non-oscillatory (WENO) schemes. Li et al. [12] developed a shifting mesh algorithm specifically to allow for the investigation of long-time rivulet formation. For an extensive review of numerical methods, we refer to the introduction section of [12].

It must be noticed that fully resolving the advancing thin film front requires a dense mesh and heavy computational costs for long-time simulations. Therefore, it is beneficial to utilize adaptive mesh method to improve the numerical accuracy and efficiency. In previous works, the $h$-adaptive (local refinement) method and $r$-adaptive (moving mesh refinement) method were commonly used in thin film flow areas. Sun et al. [13] developed an $h$-adaptive mesh refinement method based on an optimal interpolation error estimate for a 2D thin film equation in the mixed finite element formulation. Li et al. [14] developed an $h$-adaptive finite difference essentially non-oscillatory scheme for a nonlinear time-dependent gravity-driven thin film equation. It was shown that the adaptive multigrid offers increased flexibility together with a significant reduction in memory requirement.

The moving mesh method for thin film flow problem has been studied recently. Alharbi and Naire [15] worked out the moving mesh finite difference method by coupling the moving mesh partial differential equations (MMPDEs) [16] with the 1D thin film flow equations which model a surfactant-laden drop. By adapting the curvature monitor function to include multiple solution components, the moving mesh method accurately resolves the complicated multiple wave-like structures in both variables: the film height and surfactant concentration. In [17], Alharbi further studied the 2D gravity-driven fingering instability of liquid sheet by utilizing the MMPDE4 [16] and the parabolic Monge-Ampére (PMA) equation which is based on optimal transport. The results showed that the moving mesh methods are accurate and offer significant reductions in memory requirements.

---

*Corresponding author.
E-mail addresses: H.Zhang4@uu.nl (H. Zhang), P.A.Zegeling@uu.nl (P. A. Zegeling)

Owing to the great advantage of adaptive methods over fixed mesh methods, the objective of this work is to study the fourth-order thin film flow problem by using the mixed finite element formulation on the adaptive moving mesh. The moving mesh strategy used in this work is originally proposed in [18, 19] and has been widely extended and applied to many applications. For example, it has been successfully utilized in the 1D and 2D hyperbolic systems of conservation laws [20], the spike dynamics of the singularity perturbed Gierer-Meinhardt model in 2D [21], the dendritic growth in 2D and 3D governed by a phase-field model [22] and the simulations of Kohn-Sham equation [23]. A review of this moving mesh strategy may be found in [24].

The remainder of this paper is organized as follows. Section 2 introduces the thin film flow equation and the moving reference framework. Sections 3 presents the mixed finite element formulation and discusses the discretization in the space direction and the implicit-explicit scheme in the time direction. In Section 4, we present the moving mesh strategy, the choice of monitor function and the smoothing mechanism. In section 5, several numerical experiments are carried out to demonstrate the effectiveness and efficiency of the proposed scheme. Finally, Section 6 ends with some conclusions.

## 2. Mathematical model

In this section we present a brief derivation of the thin film flow equation. For more details of the background and the derivation we refer the interested readers to [1, 7, 8].

The thin film flow equation can be derived using a lubrication or long wavelength approximation of the Navier-Stokes (NS) equations [8]. Consider the flow of a fluid film down an inclined plane, let $u$ denote the thickness of the film, $x$ be the coordinate orthogonal to the gravity direction in the plane and $z$ be the coordinate down the gradient of the inclined plane. The mass conservation law in a two-dimensional situation reads

$$\frac{\partial u}{\partial t} + \nabla \cdot (u\vec{v}) = 0. \tag{1}$$

In the lubrication approximation [25], the velocity $\vec{v}$ averaging over the thickness of the film, is given by

$$\vec{v} = -\frac{u^2}{3\mu}(\nabla p - \rho g \sin\theta \vec{e}_z), \tag{2}$$

where $\mu$ is the dynamic viscosity of the fluid, $p$ is the pressure, $\rho$ is the density of the fluid, $g$ is the gravitational constant, $\theta$ is the angle of inclination (from the horizontal) and $\vec{e}_z$ is the unit vector in the $z$ direction. The Laplace-Young boundary condition [26] describes that, at the fluid interface, the pressure can be given as the difference between the component of gravity normal to the incline and the surface tension,

$$p = \rho g u \cos\theta - \gamma_0 \Delta u, \tag{3}$$

where $\gamma_0$ is the surface tension coefficient and $\Delta u$ is an approximation to the surface curvature. By substituting (2) and (3) into (1) we get a fourth-order nonlinear parabolic partial differential equation (PDE)

$$\frac{\partial u}{\partial t} + \nabla \cdot \Big[\frac{\rho g u^3}{3\mu}(\sin\theta \vec{e}_z - \cos\theta \nabla u + \frac{\gamma_0}{\rho g}\nabla(\nabla^2 u))\Big] = 0. \tag{4}$$

After a suitable rescaling [8], Eq. (4) can be written in a simplified form as

$$\frac{\partial u}{\partial t} + \frac{\partial F(u)}{\partial z} - \beta \nabla \cdot [K(u)\nabla u] + \gamma \nabla[K(u)\nabla \Delta u] = 0. \tag{5}$$

where $F(u) = K(u) = u^3$, $\beta$ is related to the vertical component of gravity and $\gamma$ measures the surface tension.

### 2.1. Moving framework

In numerical simulations, the region of interest is near the front of the thin film where effects like bumps occur. With a fixed reference frame, the spatial domain would need to be taken as the entire domain over which the flow would evolve, leading to large portions of the area where no change occurs. Therefore, we will address this issue using a moving reference framework by adding an extra term $-s\frac{\partial u}{\partial z}$ to Eq. (5), where the speed of the moving framework is the traveling wave speed given by the Rankine-Hugoniot condition using the boundary values,

$$s = \frac{k_r(u_+) - k_r(u_-)}{u_+ - u_-}. \tag{6}$$

Adding the term $-s\frac{\partial u}{\partial z}$ to Eq. (5) gives

$$\frac{\partial u}{\partial t} + \frac{\partial \hat{F}(u)}{\partial z} + \beta \nabla \cdot [K(u)\nabla u] + \gamma \nabla \cdot [K(u)\nabla \Delta u] = 0, \tag{7}$$

where $\hat{F}(u) = F(u) - su$.

The effects of the moving framework for fourth-order equations have been studied by [11, 27], their results show that the moving framework gives consistent solutions as the fixed framework.

## 3. The mixed finite element Discretization

There are several different finite element methods to solve fourth-order PDEs. One is the standard finite element method (FEM) with high order piecewise polynomials which includes derivatives as degrees of freedom. Another approach is to use a low order basis such as the mixed finite element method [13].

As in Ref. [13], we apply a mixed finite element method to the fourth-order equation. We describe the main idea briefly for the ease of reference. By defining the potential variable $w = \beta u - \gamma \Delta u$, Eq. (7) can be split into a set of second order equations,

$$\begin{cases} \frac{\partial u}{\partial t} - \nabla \cdot [K(u)\nabla w] + \frac{\partial \hat{F}(u)}{\partial z} = 0, & (x,z) \in \Omega, \\ -\beta u + \gamma \Delta u + w = 0, & (x,z) \in \Omega, \end{cases} \tag{8}$$

with initial and boundary conditions

$$\begin{cases} u(x,z;t=0) = u_0(x,z), & (x,z) \in \Omega, \\ u(x,z;t) = u_B(x,z;t), & (x,z) \in \partial\Omega_D, \\ \vec{n} \cdot \nabla u = 0, & (x,z) \in \partial\Omega/\partial\Omega_D, \\ \vec{n} \cdot \nabla w = 0, & (x,z) \in \partial\Omega, \end{cases} \tag{9}$$

where $\Omega \in \mathrm{R}^2$ is the physical domain, $\Omega_D$ is the Dirichlet boundary, $\vec{n}$ denotes outward normal direction of $\partial\Omega$.

The weak form for the mixed formulation (8) is to find $u, w \in H^1(\Omega)$ such that

$$\begin{cases} (\frac{\partial u}{\partial t}, \phi) + K(u)(\nabla w, \nabla \phi) - (\hat{F}(u), \frac{\partial \phi}{\partial z}) = 0, \\ -\beta(u, \phi) - \gamma(\nabla u, \nabla \phi) + (w, \phi) = 0, \end{cases} \tag{10}$$

where $(\cdot, \cdot)$ denotes the $L^2(\Omega)$ inner product and $\phi \in H^1(\Omega)$ is the test function for both equations.

Now we discretize the system (8) in space using the standard piece-wise linear finite element method (FEM) and we use the same space for both variables.

Denote $\mathcal{T}$ as the triangular mesh for the physical domain $\Omega$ and $h$ be the characteristic length of the triangle edge. The finite element space $P(\mathcal{T}) \in H^1(\Omega)$ is chosen as a standard linear finite element space. Let $\vec{x} = \{\vec{x}_i\}_{i=1}^N$ be the nodes of $\mathcal{T}$, the finite element space $P(\mathcal{T})$ is expanded by the basis functions $\{\phi_j\}_{j=1}^N$ such that $\phi_j(\vec{x}_i) = \delta_{ij}$, where $\delta_{ij}$ is the Kronecker delta operator. Then the approximations to $u$ and $w$ can be represented in the following form,

$$u_h = \sum_{j=1}^N u_j(t)\phi_j(x,y), \quad w_h = \sum_{j=1}^N w_j(t)\phi_j(x,y). \tag{11}$$

The time integration of the thin film flow equation is demanding because of the appearance of nonlinear fourth-order term. For stability, an explicit scheme requires a time step $\Delta t$ of the order $h^4$. The can make the explicit schemes prohibitively expensive. Another common approach for solving equation of this type is to use the first order accurate semi-implicit scheme. In this work, we apply the well known implicit-explicit (IMEX) scheme [28] to the thin film flow equation. We divide the time interval $[0, T]$ in to $N_t$ intervals of size $\Delta t > 0$ with $T = N_t \cdot \Delta t$. By treating the diffusion and hyper-diffusion terms implicitly and the nonlinear convection term explicitly, the IMEX scheme reads

$$\begin{cases} (u_h^{n+1} - u_h^n, \phi_i) + \Delta t K(u_h^{n+1})(\nabla w_h^{n+1}, \nabla \phi_i) - \Delta t(\hat{F}(u_h^n), \frac{\partial \phi_i}{\partial z}) = 0, & i = 1, 2, \cdots N, \\ -\beta(u_h^{n+1}, \phi_i) - \gamma(\nabla u_h^{n+1}, \nabla \phi_i) + (w_h^{n+1}, \phi_i) = 0, & i = 1, 2, \cdots, N. \end{cases} \quad (12)$$

Let $\bar{u} = [u_1, u_2, \cdots, u_N]^T$, $\bar{w} = [w_1, w_2, \cdots, w_N]^T$, $\bar{0}$ denote the $2N$-dimensional zero vector and $\bar{x} = [\bar{u}; \bar{w}]$ denote the combined unknown coefficients, by substituting (11) into (12), the discretized system can be presented in the following form

$$\bar{f}(\bar{x}^{n+1}) := \begin{bmatrix} \bar{f}_0(\bar{u}^{n+1}) \\ \bar{f}_1(\bar{w}^{n+1}) \end{bmatrix} = \begin{bmatrix} \mathbf{M}_1 \bar{u}^{n+1} - \Delta t \bar{F}(\bar{u}^n) + \Delta t \mathbf{K}_1(u_h^{n+1}) \bar{u}^{n+1} - \mathbf{M}_1 \bar{u}^n \\ -\beta \mathbf{M}_1 \bar{u}^{n+1} - \gamma \mathbf{K}_2 \bar{u}^{n+1} + \mathbf{M}_1 \bar{w}^{n+1} \end{bmatrix} = \bar{0}, \quad (13)$$

where the elements of $\mathbf{M}_1$, $\bar{F}(\bar{u}^n)$, $\mathbf{K}_1(u_h^{n+1})$ and $\mathbf{K}_2$ are given by

$$\mathbf{M}_{1ij} = \int_{\mathcal{T}} \phi_i \phi_j \mathrm{d}x \mathrm{d}z, \qquad \bar{F}(\bar{u}^n)_i = \int_{\mathcal{T}} \hat{F}(u_h^n) \frac{\partial \phi_i}{\partial z} \mathrm{d}x \mathrm{d}z,$$

$$\mathbf{K}_{1ij}(u_h^{n+1}) = \int_{\mathcal{T}} K(u_h^{n+1})(\nabla \phi_i, \nabla \phi_j) \mathrm{d}x \mathrm{d}z, \qquad \mathbf{K}_{2ij} = \int_{\mathcal{T}} (\nabla \phi_i, \nabla \phi_j) \mathrm{d}x \mathrm{d}z, \qquad i, j = 1, 2, \cdots, N.$$

Because of the term $K(u_h^{n+1})$ in (12), in every time step we have to solve a nonlinear system. Here we choose the well known quasi-Newton method. By evaluating $\mathbf{K}_1$ at $t^n$, an approximation to the Jacobian of (13) reads

$$\mathbf{J} = \begin{bmatrix} \mathbf{M}_1 & \Delta t \mathbf{K}_1(u_h^n) \\ -\beta \mathbf{M}_1 - \gamma \mathbf{K}_2 & \mathbf{M}_1 \end{bmatrix}. \quad (14)$$

Let the initial guess at every time step be $\bar{x}^{n+1,0} = \bar{x}^n$, the iteration of quasi-Newton method reads

$$\bar{x}^{n+1,s+1} = \bar{x}^{n+1,s} - \mathbf{J}^{-1} \bar{f}(\bar{x}^{n+1,s}). \quad (15)$$

In every iteration step the stopping criterion is chosen as $\|\bar{x}^{n+1,s+1} - \bar{x}^{n+1,s}\| \leq \text{tol}$.

*3.1. Precondition strategy*

The linear system in quasi-Newton iteration is usually very stiff because of the diffusion terms. Recently, different precondition strategies for solving fourth-order PDEs like (5) have been proposed. For example, the precondition strategies in [29, 30, 31] for Cahn-Hilliard equations.

In order to solve the linear system (15) efficiently, we apply a block Schur complement preconditioner to the block Jacobian matrix $\mathbf{J}$ and choose a generalized minimal residual (GMRES) solver for the preconditioned linear system in the quasi-Newton iterations. In the Schur complement preconditioning, we apply the approximation strategy proposed in [30].

The block-triangular preconditioner for $\mathbf{J}$ is

$$\mathbf{P} = \begin{bmatrix} \mathbf{M}_1 & 0 \\ -\beta \mathbf{M}_1 - \gamma \mathbf{K}_2 & \mathbf{S} \end{bmatrix}, \quad (16)$$

where the Schur complement of $\mathbf{J}$ is $\mathbf{S} = \mathbf{M}_1 + (\beta \mathbf{M}_1 + \gamma \mathbf{K}_2)\mathbf{M}_1^{-1}(\Delta t \mathbf{K}_1)$. In the preconditioning, in order to reduce the computational cost, we adopt the strategy proposed in [30] to approximate the Schur complement,

$$\hat{\mathbf{S}} = \hat{\mathbf{S}}_1 \mathbf{M}_1^{-1} \hat{\mathbf{S}}_2 = (\mathbf{M}_1 + \sqrt{\Delta t}(\beta \mathbf{M}_1 + \gamma \mathbf{K}_2))\mathbf{M}_1^{-1}(\mathbf{M}_1 + \sqrt{\Delta t} \mathbf{K}_1), \quad (17)$$

the approximated preconditioner is then denoted by $\hat{\mathbf{P}}$. As suggested by [30], the algebraic multigrid (AMG) preconditioner is chosen for the approximation to the inverse of $\hat{\mathbf{S}}_1$ and $\hat{\mathbf{S}}_2$.

Now, we briefly illustrate the implementation of preconditioning. The preconditioner is built to operate on the Jacobian matrix $\mathbf{J}$ in block matrix form, such that the product matrix

$$\hat{\mathbf{P}}^{-1} \mathbf{J} = \begin{bmatrix} \mathbf{M}_1^{-1} & 0 \\ \hat{\mathbf{S}}^{-1}(\beta \mathbf{M}_1 + \gamma \mathbf{K}_2)\mathbf{M}_1^{-1} & \hat{\mathbf{S}}^{-1} \end{bmatrix} \begin{bmatrix} \mathbf{M}_1 & \Delta t \mathbf{K}_1 \\ -\beta \mathbf{M}_1 - \gamma \mathbf{K}_2 & \mathbf{M}_1 \end{bmatrix}, \quad (18)$$

is of a from that Krylov subspace-based iterative solver like GMRES can solve in a few iterations.

In preconditioning, we need to solve $[\bar{x}_0; \bar{x}_1] = \hat{\mathbf{P}}^{-1}[\bar{f}_0; \bar{f}_1]$. We first get $\bar{x}_0$ using $\bar{x}_0 = \mathbf{M}_1^{-1} \bar{f}_0$ and then compute $\bar{x}_1 = \hat{\mathbf{S}}^{-1}(\bar{f}_1 + (\beta \mathbf{M}_1 + \gamma \mathbf{K}_2)\bar{x}_0)$.

## 4. Mesh redistribution and monitor smoothing

In this section we briefly introduce the moving mesh strategy proposed by [18, 19], in which the physical equation and the moving mesh equation are solved alternately. In particular, we will concentrate on the choice of monitor function smoothing method.

*4.1. Moving mesh strategy*

Let $\vec{x}$ and $\vec{\xi}$ denote the coordinates of the physical domain $\Omega$ and the computational domain $\Omega_c$, respectively. In [18], a harmonic mapping $\vec{\xi} = \vec{\xi}(\vec{x})$ from $\Omega$ to $\Omega_c$ is achieved by solving the Euler-Lagrange equation

$$\frac{\partial}{\partial x^i}\left(G^{ij}\frac{\partial \xi^k}{\partial x^j}\right) = 0, \tag{19}$$

where $M = (G^{ij})^{-1}$ is the monitor function.

In order to give reasonable distribution of grid points on the boundary, Li et. al. [19] introduced an optimization problem with some appropriate constraints to redistribute the interior and boundary points simultaneously,

$$\begin{cases} \min \quad E(\vec{\xi}) \\ \text{s.t.} \quad \vec{\xi}|_{\partial\Omega} = \vec{\xi}_b \in \mathcal{K}, \end{cases} \tag{20}$$

where $\mathcal{K}$ is an admissible set for the boundary mappings and $E(\vec{\xi})$ is the mesh energy defined by

$$E(\vec{\xi}) = \sum_k \int_\Omega G^{ij}\frac{\partial \xi^k}{\partial x^i}\frac{\partial \xi^k}{\partial x^j}\mathrm{d}\vec{x}. \tag{21}$$

Ref. [19] proposed an iterative algorithm to move the mesh. In the beginning, the initial mesh $\vec{\xi}^0$ on $\Omega_c$ is generated by solving the following optimization problem:

$$\begin{cases} \min \quad \sum_k \int_\Omega \sum_i \left(\frac{\partial \xi^k}{\partial x^i}\right)^2 \mathrm{d}\vec{x}, \\ \text{s.t.} \quad \vec{\xi}|_{\partial\Omega} = \vec{\xi}_b \in \mathcal{K}. \end{cases} \tag{22}$$

At time step $t^n$, denote the mesh generated from problem (20) as $\vec{\xi}^{(n)}$, the moving mesh algorithm can be summarized as

1. Get $\vec{\xi}^n$ by solving (20) and compute the difference $\delta\vec{\xi} = \vec{\xi}^0 - \vec{\xi}^n$. If $\|\delta\xi\|_{L^\infty}$ is smaller than a given tolerance, then the mesh-redistribution at time step $t^n$ is finished. Otherwise, do step 2 to step 4.
2. Obtain the displacement of the physical mesh $\delta\vec{x}$ by using $\delta\vec{\xi}$ and the Jacobi matrix, then move the physical mesh by

$$\vec{x}^n = \vec{x}^n + \tau\delta\vec{x}^n, \tag{23}$$

   where $\tau \in [0,1]$ is a parameter to prevent mesh tangling.
3. Update the solution $u_h^n$ on the new mesh $\vec{x}^n$.
4. Update the monitor function using the updated solution $u_h^n$ and go to step 1.

For more details of the moving mesh strategy and the solution update algorithm we refer the readers to the Refs. [18, 19].

*4.2. Choice of the monitor function*

In the moving mesh method, the monitor function $M$ connecting the mesh with the physical solution, is chosen to redistribute more gird points at critical regions where more accuracy is needed, thereby reducing errors introduced by the numerical scheme. In this work, we consider an adaptive monitor function [32, 33, 34]

$$M = (1-\kappa)\gamma(u) + \kappa\omega, \tag{24}$$

where $\gamma(u)$ is a normalization of the monitor component $\omega$,

$$\gamma(u) = \frac{1}{|\Omega|} \iint \omega \mathrm{d}x \mathrm{d}z,$$

and the parameter $\kappa$ indicate the ratio of points in the critical regions [35].

In practice, the monitor component $\omega$ can be chosen as a function of the gradient of $u$:

$$\omega = |\nabla u|, \tag{25}$$

or a function of the second order derivative of $u$:

$$\omega = |\Delta u|^{\frac{1}{2}}. \tag{26}$$

In the following, we call the monitor using (25) as the arc-length type monitor and the monitor using (26) as the curvature type monitor.

### 4.3. Smoothing mechanism

Since the computed monitor $M$ is usually not smooth, in order to avoid a very distorted mesh around critical regions, the monitor function is generally smoothed [36, 37, 38] before the solving of the moving mesh problem. In the moving mesh finite element method (MMFEM) situation, one approach to smooth the monitor is to filter the monitor several times [18], but this approach needs to determine the filter time and may be too costly if filter time is big. Instead, we apply a smoothing strategy based on a diffusive mechanism in [37]. Similar smoothing strategies have also been adopted in [22, 39, 40, 41] and obtained good results.

A 2D extension of the smoothing equation in [37] is given by

$$\begin{cases} [\mathcal{I} - (\sigma_\xi(\sigma_\xi + 1)(\Delta\xi)^2 \frac{\partial^2}{\partial \xi^2} + \sigma_\eta(\sigma_\eta + 1)(\Delta\eta)^2 \frac{\partial^2}{\partial \eta^2})]\tilde{M} = M, & (\xi, \eta) \in \Omega_c, \\ \vec{n} \cdot \nabla \tilde{M} = 0, & (\xi, \eta) \in \partial\Omega, \end{cases} \tag{27}$$

where $\mathcal{I}$ is the identity operator, $\sigma_\xi$ and $\sigma_\eta$ are the spatial smoothing parameters in $\xi$- and $\eta$-directions. This smoothing equation is defined on the initial computational mesh which is fixed for all time steps, therefore, we only need to calculate the discretization of the linear operator once and used it for all later steps. At every time step, we preform several conjugate gradient (CG) iterations to obtain reasonable approximation to $\tilde{M}$. With this smoother monitor $\tilde{M}$, the corresponding mesh will be less singular, hence the physical equation can be solved more efficiently.

### 4.4. Properties of the adaptive mesh

Before we solve the smoothing equation, we would like to study the effects of the spatial smoothing parameters $\sigma_\xi$ and $\sigma_\eta$ in (27). Refs. [37, 42] studied a 1D moving mesh partial differential equation (MMPDE),

$$\begin{cases} \frac{\partial}{\partial \xi}\left(\frac{\dot{\tilde{n}}}{\omega}\right) = -\frac{1}{\tau_s}\frac{\partial}{\partial \xi}\left(\frac{\tilde{n}}{\omega}\right), \\ \tilde{n} = [\mathcal{I} - \sigma_s(\sigma_s + 1)(\Delta\xi)^2 \frac{\partial^2}{\partial \xi^2}]n, \end{cases} \tag{28}$$

the authors proved that the 1D mesh obtained using this MMPDE admits a local quasi-uniformity: $|\frac{x_{\xi\xi}}{x_\xi}| \leq \frac{1}{\sqrt{\sigma(\sigma+1)}\Delta\xi}$ with discretized version:

$$\frac{\sigma}{\sigma + 1} \leq \frac{\Delta x_{i+1}(t)}{\Delta x_i(t)} \leq \frac{\sigma + 1}{\sigma}, \quad \forall t \in [0, T]. \tag{29}$$

Although we are not able to prove similar properties for the adaptive mesh obtained using the moving mesh strategy with monitor smoothed by the 2D extension (27), we numerically show that the adaptive mesh has similar properties:

$$\begin{cases} \text{Density ratio along } x\text{-direction:} & \frac{\sigma_\xi}{\sigma_\xi + 1} \leq \frac{\Delta x_{i+1}(t)}{\Delta x_i(t)} \leq \frac{\sigma_\xi + 1}{\sigma_\xi}, \\ \text{Density ratio along } z\text{-direction:} & \frac{\sigma_\eta}{\sigma_\eta + 1} \leq \frac{\Delta z_{j+1}(t)}{\Delta z_j(t)} \leq \frac{\sigma_\eta + 1}{\sigma_\eta}, \quad \forall t \in [0, T]. \end{cases} \tag{30}$$

Since we are working on 2D triangular mesh, we approximate $\frac{\Delta x_{i+1}(t)}{\Delta x_i(t)}$ ($\frac{\Delta z_{j+1}(t)}{\Delta z_j(t)}$) using the area ratio $\frac{\text{area}_{i+1}}{\text{area}_i}$ along $x$-direction ($\frac{\text{area}_{j+1}}{\text{area}_j}$ along $z$-direction).

Choosing the test function as

$$u(x, z) = -\tanh(100z)\tanh(100x), \quad (x, z) \in [-0.5, 0.5] \times [-0.5, 0.5]. \tag{31}$$

The surface plot in Fig. 1 indicates that the gradient norm of this function is higher along the $x$ and $z$ axes than in other regions.

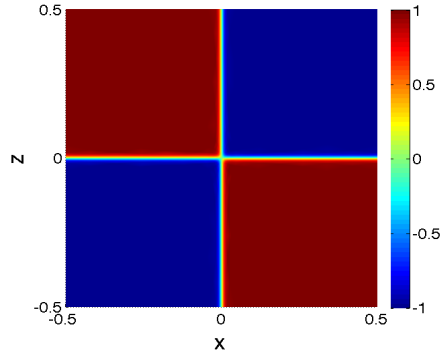

Figure 1: Surface plot of the test function (31).

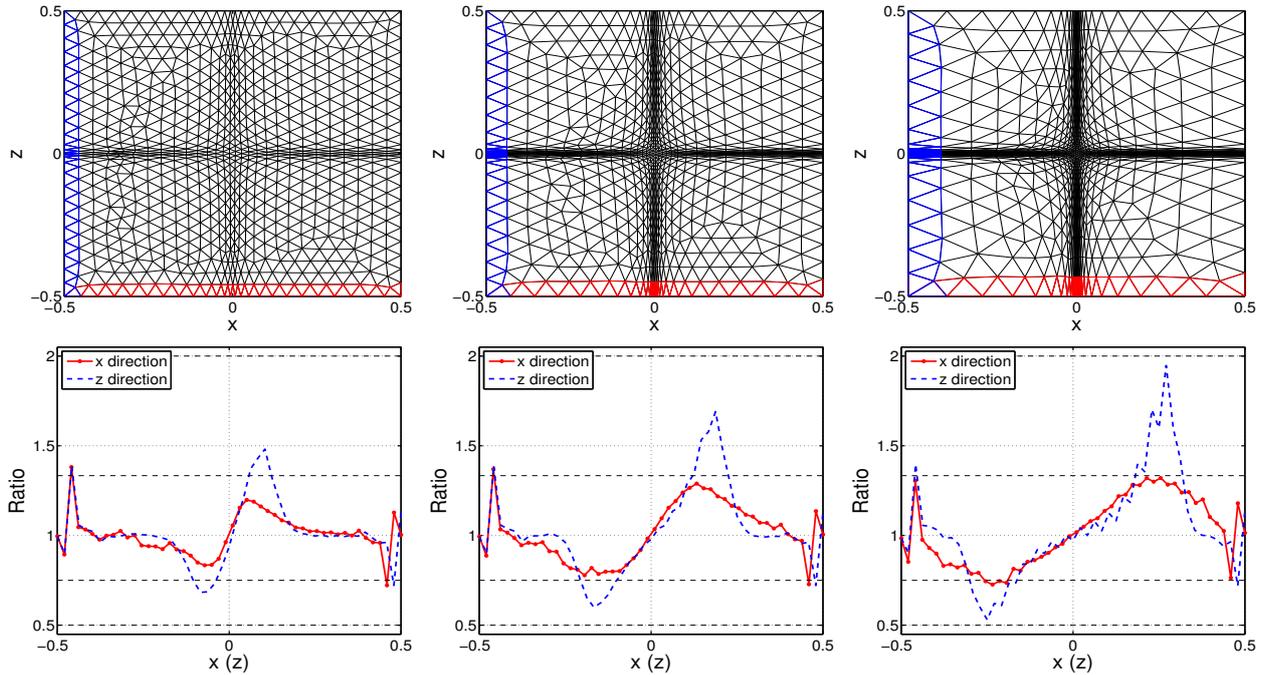

Figure 2: Adaptive meshes and density ratios for $\kappa = 0.25$ (left), $\kappa = 0.5$ (middle) and $\kappa = 0.75$ (right) with $\sigma_\xi = 3, \sigma_\eta = 1$. Red dotted and blue dashed curves denote the density ratios along the horizontal and vertical boundaries, respectively. The dashed lines from bottom to top denote density ratios $\frac{1}{2}, \frac{3}{4}, \frac{4}{3}, 2$.

We use the arc-length type monitor and fix the spatial smoothing parameters $\sigma_\xi = 3, \sigma_\eta = 1$. By choosing different adaptivity parameter $\kappa = 0.25, 0.5, 0.75$, the corresponding meshes and density ratios at the horizontal and vertical boundaries of the meshes are present in Fig. 2. We can observe that the meshes have different density ratios in different directions. For all values of $\kappa$, the density ratios of the horizontal and vertical boundaries are bounded by the interval $[\frac{\sigma_\xi}{\sigma_\xi+1}, \frac{\sigma_\xi+1}{\sigma_\xi}]$ and $[\frac{\sigma_\eta}{\sigma_\eta+1}, \frac{\sigma_\eta+1}{\sigma_\eta}]$ (the violations at the corners are due to the non-uniformity of the initial triangular mesh), respectively. With the increase of adaptivity parameter $\kappa$, the minimum and maximum density ratios get more and more close to the lower and upper bounds. These plots clearly show the local quasi-uniformity of the adaptive mesh obtained by the smoothed monitor.

Table 1: Space accuracy test of the MMFEM using the fixed mesh and the moving mesh ($T = 0.01$, $\Delta t = 1.0\text{e-}5$).

| Mesh size | Fixed mesh | | Moving mesh | | | |
| --- | --- | --- | --- | --- | --- | --- |
| | | | arc-length type | | curvature type | |
| | $L^2$ error | $L^2$ order | $L^2$ error | $L^2$ order | $L^2$ error | $L^2$ order |
| $h = 0.1$ | 1.0731e-02 | − | 1.1289e-02 | − | 1.1370e-02 | − |
| $h = 0.05$ | 2.5489e-03 | 2.0743 | 2.9667e-03 | 1.9280 | 2.8239e-03 | 2.0095 |
| $h = 0.025$ | 5.9273e-04 | 2.1040 | 7.5613e-04 | 1.9722 | 6.8194e-04 | 2.0500 |
| $h = 0.0125$ | 1.7112e-04 | 1.7924 | 1.6521e-04 | 2.1943 | 1.4097e-04 | 2.2743 |

## 5. Numerical experiments

In this section we present some numerical results computed using the MMFEM described in the previous sections. First, we show the numerical convergence order of the scheme. Then we show the accuracy, efficiency and different features of the moving mesh method. At last, we present the simulations of 2D finger pattern.

Our codes are based on the AFEPack [43] and we use the 2D mesh generator EasyMesh [44] to generate triangular meshes. In order to implement the AMG preconditioning, we also used some packages from deal.II [45] which provide wrapper classes to use the linear algebra parts of the Trilinos library [46]. For all AMG preconditioners we choose two steps of Chebyshev smoother and two V-cycles. The tolerance used in the quasi-Newton iteration is take to be tol = 1.0e-6 and the tolerance used in the GMRES solver is 1.0e-8 . For all numerical simulations we choose $\sigma_\xi = \sigma_\eta = 1$.

### 5.1. Convergence order

To test the convergence order of the scheme (12) in the space direction, we consider a linear equation whose solution is smooth:

$$\begin{cases} \frac{\partial u}{\partial t} - \beta \Delta u + \gamma \Delta^2 u = 0, & (x,y) \in \Omega = [0,1] \times [0,1], \\ u(x,z;t=0) = \cos(2\pi x)\cos(2\pi z), \\ \vec{n} \cdot \nabla u = 0, & (x,z) \in \partial\Omega. \end{cases} \quad (32)$$

The exact solution of (32) is

$$u_{ex}(x,z;t) = \exp(-(8\beta\pi^2 + 64\gamma\pi^4)t)\cos(2\pi x)\cos(2\pi z).$$

Choosing $\beta = 0.5, \gamma = 0.0025$, $\sigma_\xi = \sigma_\eta = 1$, $\kappa = 0.5$ and fixing the time step size $\Delta t = 1.0\text{e-}5$, the corresponding errors and orders for $h = 0.1, 0.05, 0.025, 0.0125$ at $T = 0.01$ are presented in Table 1, from which a second-order rate of convergence is observed.

In Fig. (3) we present the adaptive meshes obtained using both the acr-length type monitor and the curvature type monitor with initial mesh size $h = 0.0125$. Because of the smoothness of the cosine shape solution, the effect of the moving mesh method is weak, but we can still observe some adaptiveness from the meshes.

### 5.2. Simulations of traveling waves

Analytical solutions to the nonlinear thin film flow equation are not always available, so we will use the MMFEM to produce one-dimensional traveling wave (TW) solutions to Eq. (5) and verify the performance of the method. Since the flux only exists in the $z$-direction, we consider the TW solutions in this direction. The TW connecting two regions of different values has the form

$$\begin{cases} u(z;t) = u(\zeta) = u(z - st), \\ u(-\infty) = u_-, \quad u(+\infty) = u_+. \end{cases}$$

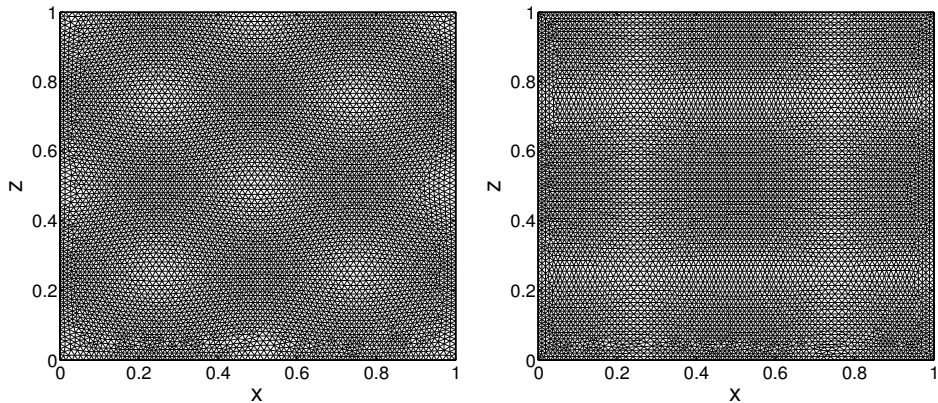

Figure 3: Adaptive meshes for (32) obtained using the arc-length type monitor (left) and the curvature type monitor (right) with initial mesh size $h = 0.0125$.

Assuming all the derivatives of the TW vanish as $\zeta \to \pm\infty$, by substituting $u(\zeta)$ into (5) we obtain a fourth-order ordinary differential equation (ODE)

$$-su' + F'(u) - \beta[K(u)u']' + \gamma[K(u)u''']' = 0. \tag{33}$$

Integrating (33) once gives

$$-s(u - u_+) + F(u) - F(u_+) - \beta K(u)u' + \gamma K(u)u''' = 0, \tag{34}$$

subject to the boundary conditions

$$u(-\infty) = u_-, \quad u(+\infty) = u_+. \tag{35}$$

The equation (34) can be treated as either an initial value problem (IVP) or as a boundary value problem (BVP). In this work we solve (34) as a BVP. Since (34) is of third-order, we will impose the third boundary condition as

$$u'''(-\infty) = 0. \tag{36}$$

By introducing $v = u'$ and $w = v'$, we transform (34) to a first-order system of ODEs

$$\begin{cases} u' = v, \\ v' = w, \\ w' = \dfrac{1}{\gamma K(u)}[s(u-u_+) - (F(u) - F(u_+)) + \beta K(u)v], \\ u(\zeta_-) = u_-, \quad u(\zeta_+) = u_+, \quad w(\zeta_-) = 0, \quad \zeta \in [\zeta_-, \zeta_+]. \end{cases} \tag{37}$$

The above system is then solved using the built-in BVP solver 'bvp4c' in Matlab [47].

To study the TW solutions to the thin film flow model, we consider a specific situation [3], in which the film flows up an inclined plane because of the Marangoni stress [48] created by a temperature gradient on the planes. In this case, the counteracting effect of the Marangoni stress and the gravitational force leads to a non-convex flux function

$$F(u) = u^2 - u^3. \tag{38}$$

Taking $K(u) = u^3$, $\beta = 0$, $\gamma = 1$, the governing equation reads

$$u_t + \frac{\partial(u^2 - u^3)}{\partial z} + \nabla \cdot (u^3 \nabla \Delta u) = 0. \tag{39}$$

The significance of the non-convex flux function is that it allows (39) admitting different types of TW solutions. Bertozzi et al. [5, 6, 49] showed that (39) may have the admissible Lax shock, undercompressive shock and rarefaction wave, depending on the choice of initial conditions. In the following, we will test the MMFEM with different initial conditions and see if we can get the same numerical results as in [5, 11, 50].

Table 2: Comparison of CPU times between the fixed mesh and the moving mesh $((x,z) \in [0, 0.5] \times [0, 5], T = 20, \Delta t = 0.1)$.

| Fixed mesh | | Moving mesh ($\kappa = 0.3$) | |
|---|---|---|---|
| Mesh size | CPU time [s] | Mesh size | CPU time [s] |
| $h = \frac{1}{16}$ | 37.58 | $h = \frac{1}{8}$ | 24.80 |
| $h = \frac{1}{32}$ | 141.13 | $h = \frac{1}{16}$ | 72.09 |

In previous studies [5, 11, 50], the value of $\gamma$ is taken to be 1 and the length of the computational interval in the $z$-direction ranges from 300 to 1000 and the end time ranges from 500 to 100000. In order to simplify the computation, we introduce $\hat{x} = 0.1x, \hat{z} = 0.1z, \hat{t} = 0.1t$, then (39) can be rescaled to

$$u_{\hat{t}} + \frac{\partial (u^2 - u^3)}{\partial \hat{z}} + 0.001 \hat{\nabla} \cdot (u^3 \hat{\nabla} \hat{\Delta} u) = 0, \tag{40}$$

where $\hat{\nabla} = [\nabla_{\hat{x}}, \nabla_{\hat{z}}]^T$, $\hat{\Delta} = \frac{\partial^2}{\partial \hat{x}^2} + \frac{\partial^2}{\partial \hat{z}^2}$. Now in the moving framework, the computational interval in $\hat{z}$-direction can be fixed as, for example, $[0, 5]$, and the end time $T = 100$ is large enough to allow the initial profile to develop into a TW solution. For convenience, in the following simulations, we still use $x, z, t$ instead of $\hat{x}, \hat{z}, \hat{t}$.

Consider $u_- = 0.3323$ and $u_+ = 0.1$, we solve (40) with three different initial conditions [5] connecting $u_-$ and $u_+$:

$$\text{Case 1}: \quad u(x) = \frac{1}{2}(u_- - u_+)(1 - \tanh(10(z - 2.5))) + u_+, \tag{41}$$

$$\text{Case 2}: \begin{cases} u(x) = \frac{1}{2}(u_- - 0.6)(1 - \tanh(10(z - 2)) + 0.6, & x \leq 2.5, \\ u(x) = \frac{1}{2}(0.6 - u_+)(1 - \tanh(10(z - 3)) + u_+, & x > 2.5, \end{cases} \tag{42}$$

$$\text{Case 3}: \begin{cases} u(x) = \frac{1}{2}(u_- - 0.6)(1 - \tanh(10(z - 1.5)) + 0.6, & x \leq 2.5, \\ u(x) = \frac{1}{2}(0.6 - u_+)(1 - \tanh(10(z - 3.5)) + u_+, & x > 2.5, \end{cases} \tag{43}$$

In the simulations, we solve (40) using the moving framework with speed $s = \frac{F(u_-) - F(u_+)}{u_- - u_+} = 0.2786$ in the 2D domain $[0, 0.5] \times [0, 5]$ with Dirichlet boundary conditions on $z = 0, 5$ and Newman boundary conditions on $x = 0, 0.5$. The 1D solutions are extracted from the 2D solutions.

**Case 1**

In Fig. 4 we present the computed solutions for the fixed mesh and the moving mesh for Case 1. We choose $h = \frac{1}{16}, \frac{1}{32}$ for the fixed mesh and $h = \frac{1}{8}, \frac{1}{16}$ for the moving mesh. The time step size is taken as $\Delta t = 0.1$. Since in all three cases, $u$ is bounded by $0 \leq u \leq 2/3$, the max wave speed is $F'(1/3) = 1/3$, then the CFL numbers are in the interval $[0.2667, 1.0667]$. By taking the TW solution computed from the BVP (37) as the most accurate one, we can observe that when the meshes are refined, both the solutions and phase planes of the fixed mesh and moving mesh converge to the TW profiles. It is worth mentioning that, with the moving mesh, the solutions computed using $h = \frac{1}{8}$ is comparable to the fixed mesh solution with $h = \frac{1}{16}$. When $h = \frac{1}{16}$, the moving mesh solution almost coincides with the TW solution. The plotted grids also show that in the moving mesh situation, the arc-length type monitor redistributes more grid points along the steep front while the curvature type monitor function clusters more grid points near the overshoot and downshoot areas. Therefore, the resolutions at the overshoot and downshoot obtained by curvature type monitor are more accurate than those obtained by the arc-length type monitor. The linear stability analysis in [51, 52] states that the presence of the bump is a necessary condition for the instability of the fluid to small perturbations in the transverse direction. Therefore, in the following examples we will only apply the curvature type monitor function which gives higher resolution at the bump.

Fig. 5 shows the numerical solution (top) and the adaptive mesh (bottom) obtained by the curvature type monitor. It clearly shows that more grid points are distributed around the overshoot and downshoot areas with higher curvature values than other regions.

Table 2 gives a comparison of the CPU times used by the fixed mesh and the moving mesh with curvature type monitor. As can be seen, because of the additional cost of mesh refinement, the moving mesh takes more CPU time than the fixed mesh when using the same initial mesh size $h = \frac{1}{16}$. But to achieve solutions of the same quality, the moving mesh requires less CPU time (24.80 [s] and 72.09 [s]) than the fixed mesh (37.58 [s] and 141.13 [s]). From Fig. 4 and Table 2 we can conclude that the MMFEM is more efficient for obtaining the same solution accuracy than the fixed mesh method.

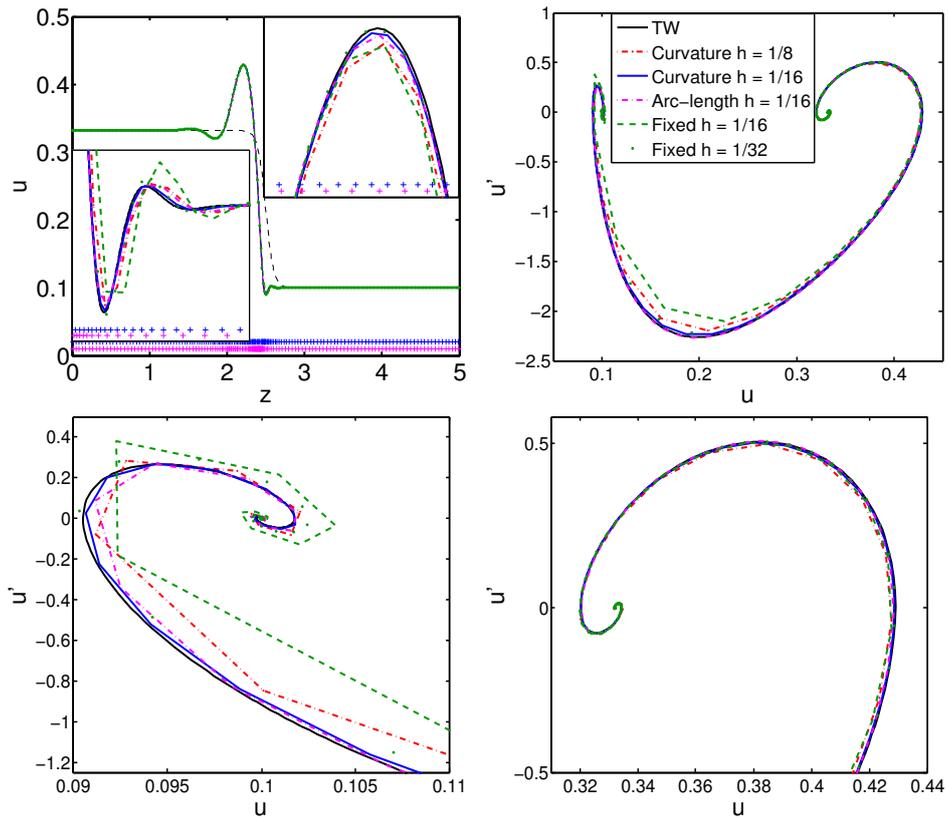

Figure 4: Case 1: Comparisions of solutions (top left), phase planes (top right), close-ups of the left part (bottom left) and right part (bottom right) of the phase planes between the fixed mesh and the moving mesh using the curvature type monitor and the arc-length type monitor ($\sigma_\xi = \sigma_\eta = 1$, $\kappa = 0.3$, $\Delta t = 0.1$, $T = 100$). The dashed curve in the top left figure is the initial condition, the cross markers at the bottom denote the 1D grid points.

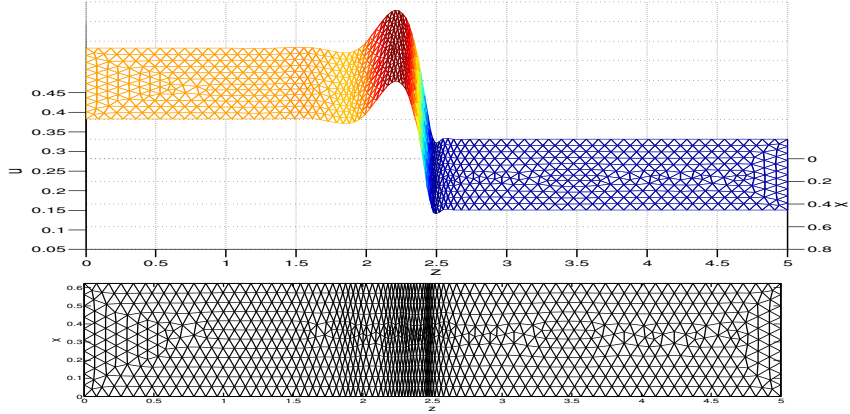

Figure 5: 2D solution (top) and adaptive mesh (bottom) at $t = 100$ obtained using the moving mesh method with curvature type monitor, $\sigma_\xi = \sigma_\eta = 1, \kappa = 0.3, \Delta t = 0.1$.

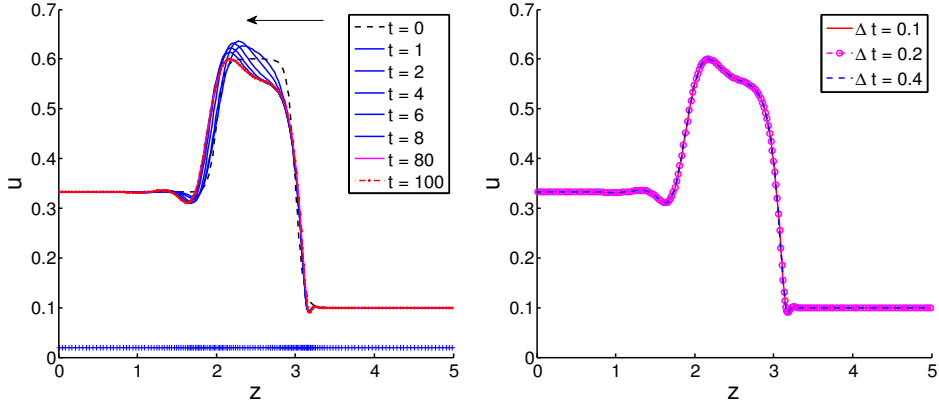

Figure 6: Solutions at different times (left) and with different time step sizes at $t = 100$ (right) using the moving mesh for case 2. Blue cross markers at the bottom of the figures denote the 1D grid points. $\sigma_\xi = \sigma_\eta = 1, \kappa = 0.3, \Delta t = 0.1$.

**Case 2**

In the second case, a bump with width 1 (corresponding to width 10 in [5]) is introduced to the initial condition. Bertozzi et al. [5] pointed out that this initial profile will develop into a different TW solution moving with the same speed as in Case 1.

In Fig. 6 we plot the evolution of the initial condition at times $t = 0, 1, 2, 4, 6, 8, 80, 100$. The initial bump evolves to a TW with an undercompressive wave on the right and a compressive wave on the left. At $t = 0$, the speeds of the left and right waves are

$$s_l = \frac{F(u_-) - F(0.6)}{u_- - 0.6} \approx 0.2625 < s, \quad s_r = \frac{F(0.6) - F(u_+)}{0.6 - u_+} = 0.2700 < s. \quad (44)$$

Both speeds are smaller than the TW wave speed $s$, so the waves move a little to the left of the initial condition. As time goes by, the speeds of both waves increase to $s$ and keep constant. Therefore, this TW solution is stable.

In Fig. 6 (right) we plot the solutions obtained using different time steps. All curves show that the TW wave is stable. This also indicates the robustness of the MMFEM.

**Case 3**

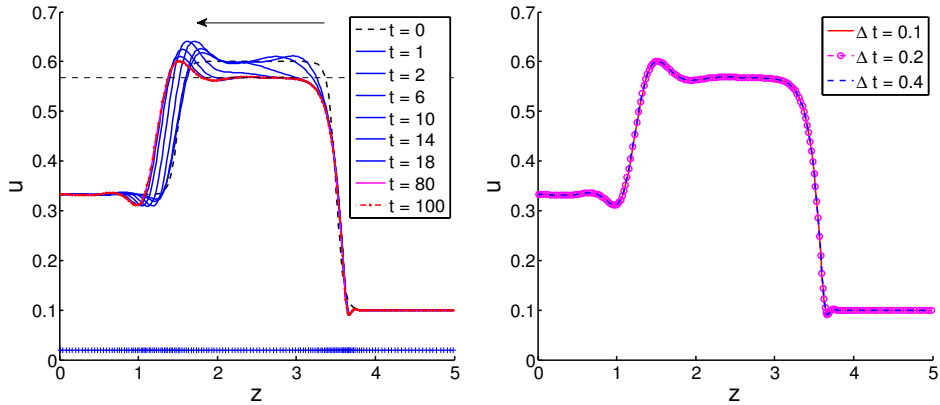

Figure 7: Solutions at different times (left) and with different time step sizes at $t = 100$ (right) using the moving mesh for Case 3. Blue cross markers at the bottom of the figures denote the 1D grid points. $\sigma_\xi = \sigma_\eta = 1, \kappa = 0.3, \Delta t = 0.1$.

In Case 3, the initial bump width is increased from 1 to 2. In [5], the authors claimed that with a larger bump width the solution does not settle down to a single traveling wave. Instead, two shocks emerge, an undercompressive wave on the right connecting $u_+$ to a larger state $u_{uc} \approx 0.568$, followed by a slower compressive wave connecting $u_{uc}$ to $u_-$. Both waves travel more slowly than the TW wave in Case 1.

Fig. 7 shows our numerical results at different time steps. As the same in Case 2, at the beginning, both waves travel a little to the left of the initial condition. As time goes by, a constant state with height $u_{uc} \approx 0.568$ appears. The speeds of the compressive wave and the undercompressive wave can then be computed as

$$s_c = \frac{F(u_-) - F(u_{uc})}{u_- - u_{uc}} \approx 0.2786, \quad s_{uc} = \frac{F(u_{uc}) - F(u_+)}{u_{uc} - u_+} \approx 0.2786. \tag{45}$$

Therefore, both waves travel at the same speed $s$, and the waves are stable rather than unstable. In [11], Ha et al. studied this case using the Crank-Nicolson scheme, the second-order Godunov scheme with limiters, the adapted upwind scheme and the WENO scheme. Among all these schemes, the WENO scheme showed very little spread of bump width for a range of CFL numbers while other schemes were sensitive to choice of the CFL number. They reached the conclusion that the TW is stable, it was merely the choice of numerical schemes and step sizes that led to the bump spreading. Fig. 7 presents solutions computed using different time steps at $T = 100$. The results again show the robustness of the MMFEM.

### 5.3. Simulations of 2D finger phenomenon

In this section we show the long time evolution of the fingering instability. Previous study [8] states that for the thin film flow equation, perturbations of long wavelengths are linearly unstable and the short wavelengths are stabilized by surface tension effects. In the simulations of finger pattern, we impose a cosinusoidal perturbation characterized by the wavelength $\lambda_0$ and the amplitude $A_0$ along the $x$-direction of the initial condition. Then we track the evolution of the 2D solution to show the evolution of finger phenomenon.

Consider the physical domain $[0, 15] \times [0, 30]$, we use the Dirichlet boundary conditions on $z = 0, 30$ and the homogeneous Neumann boundary conditions on $x = 0, 15$. The parameters in the smoothing step are $\sigma_\xi = \sigma_\eta = 1$, $\kappa = 0.5$. Let $u_- = 1$ and $u_+ = 0.1$, the perturbed initial condition is given by

$$u(x, z; 0) = 0.5(u_- - u_+)(1 - \tanh(z - 15 + A_0 \cos(2\pi x/\lambda_0))) + u_+. \tag{46}$$

Here the amplitude of the perturbation is $A_0 = 0.2$ ans the wave length $\lambda_0 = 15$. This type of perturbation method has also been adopted in [27] to simulate the finger formations in particle-laden thin film flow.

In Eq. (5) we take $F(u) = K(u) = u^3$, $\beta = 0, \gamma = 1$ and perform computations using the reference framework moving with the velocity of the unperturbed flow $s = \frac{F(u_-) - F(u_+)}{u_- - u_+}$.

In order to show the improvement of the moving mesh over the fixed mesh, we compute the solutions to the same problem with the fixed mesh (mesh size $h = 0.2, 0.4$) and the moving mesh (initial mesh size $h = 0.4$). Fig. 8 presents the back view of the thin film front at $t = 40$. By comparing these three solutions at the same time step, we can find that spurious oscillations appear at the bottom of the solution obtained by the fixed mesh with $h = 0.4$, while the moving mesh successfully redistributes more gird points near the critical regions (see Fig. 9 and the mesh in Fig. 10) and suppresses spurious oscillations.

In Fig. 10 we illustrate the surface plots of $u$ and corresponding meshes at time steps $t = 0, 20, 40, 80$ for the moving mesh method. We clearly observe the formation of a single finger. At the beginning, a bump is observed to slowly start forming. As time increases, the solution appears to develop into a single finger. The plots of $u$ also show that the speed of the finger front is higher than the TW speed $s$ while the speed at the finger root is smaller. This phenomenon agrees with the results in Ref. [7]. The obtained mesh at different time steps capture the finger shape and shows the effectiveness of the moving mesh method.

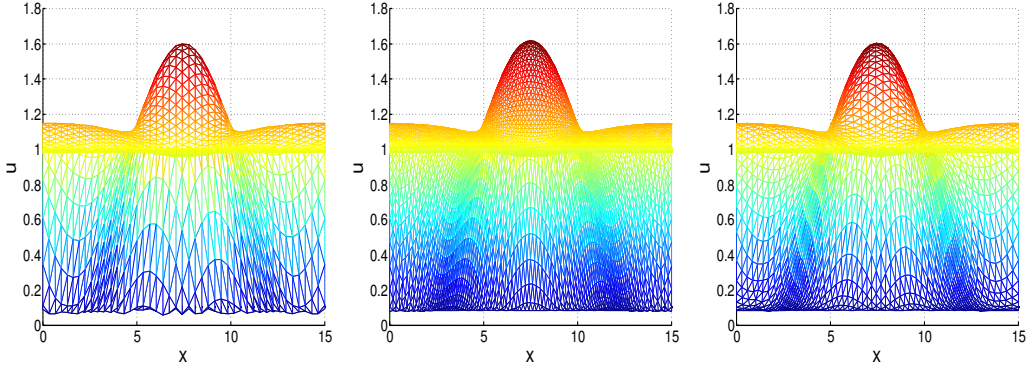

Figure 8: Back view of numerical solutions obtained by fixed meshes and moving mesh at $t = 40$, $\Delta t = 0.1$. Left: fixed mesh $h = 0.4$; middle: fixed mesh $h = 0.2$; right: moving mesh $h = 0.4$.

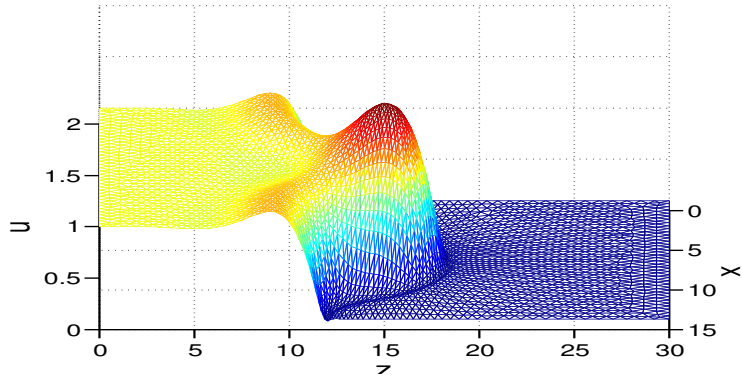

Figure 9: 3D view of the numerical solution obtained by the moving mesh method at $t = 40$, $\Delta t = 0.1$, $h = 0.4$.

## 6. Conclusions

In this paper, we have successfully solved the fourth-order thin film flow equations in a mixed finite element formulation with moving mesh refinement based on harmonic mappings [18, 19]. In order to efficiently solve the discretized nonlinear system, we used the quasi-Newton iteration method and applied the block-triangular Schur complement preconditioning with the recently developed Schur complement approximation strategy [30]. In the moving mesh step, we smoothed the monitor function using a 2D extension of the diffusive mechanism in [37]. Numerical results showed that the 2D smoothing equation admits a local quasi-uniformity which helped to reduce the singularity of the adaptive mesh. In numerical simulations, we compared the traveling wave solutions obtained by the fixed mesh and the moving mesh with the arc-length type monitor and the curvature type monitor. The results demonstrated the different features of the monitor functions: the arc-length type monitor increased the resolutions at steep regions by distributing more grid points near those regions while the curvature type monitor increased the accuracy in overshoot and downshoot regions by clustering more grid points there. The comparisons between the moving mesh and fixed mesh also showed that the moving mesh method needs shorter CPU time and less grid points to obtain solutions of the same quality than the fixed mesh method.

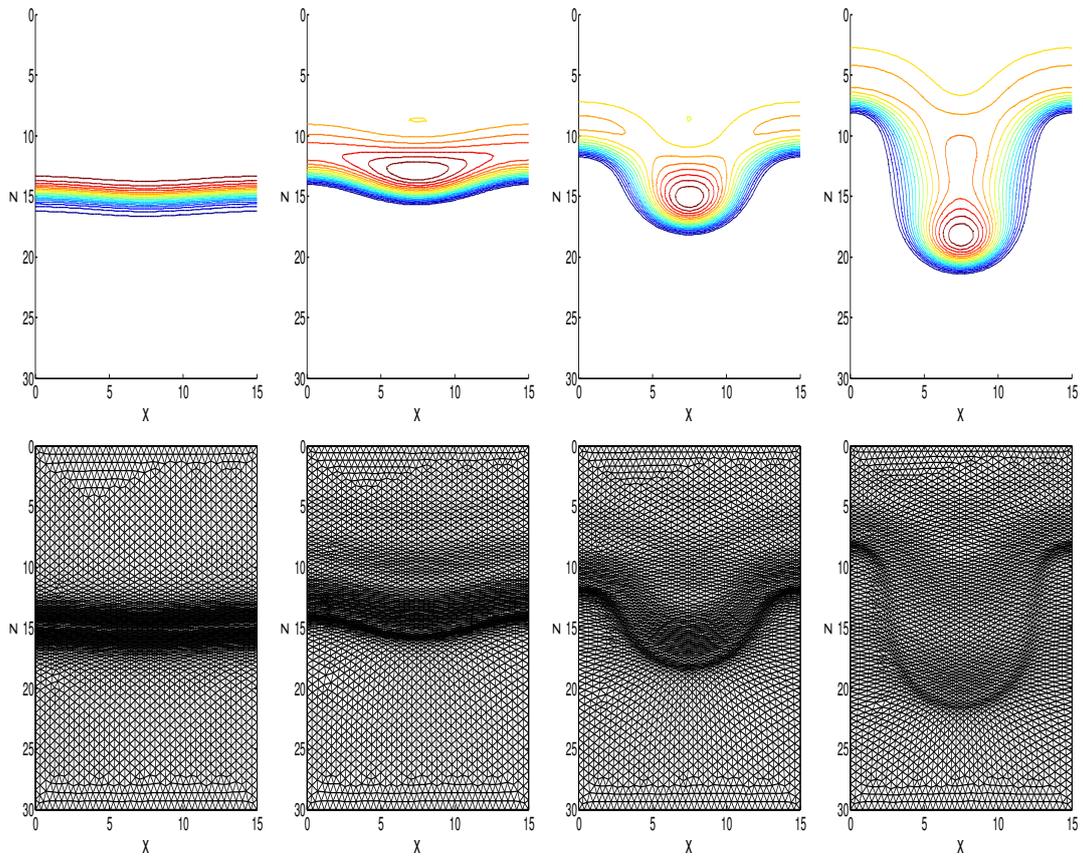

Figure 10: Numerical solutions and corresponding meshes obtained by moving mesh method with $h = 0.4$, $\Delta t = 0.1$. From left to right: $t = 0$, $t = 20$, $t = 40$, $t = 80$.

**Acknowledgements**

H. Zhang gratefully acknowledges the financial support from the China Scholarship Council (No. 201503170430). He would also like to thank Prof. Heyu Wang (Zhejiang University) for helpful discussions about solving the nonlinear algebraic system.